\documentclass[a4paper,12pt]{article}
\usepackage{amssymb}

\newcommand{\beq}{\begin{equation} }
\newcommand{\eqq}{\end{equation} }
\newcommand{\cuad}{{\sqcap\kern-.68em\sqcup}}

\newtheorem{remark}{Remark}[section]
\newcommand{\bremark}{\begin{remark} \em}
\newcommand{\eremark}{\end{remark} }

\def\beeq{\begin{equation}}
\def\eeq{\end{equation}}
\newcommand{\begeqaet}{\begin{eqnarray*}}
\newcommand{\eneqaet}{\end{eqnarray*}}

\newcommand{\foral}{\quad\mbox{for all}\quad}

\let\Section=\section
\def\section{\setcounter{equation}{0}\Section}
\newtheorem{Lem}{Lemma}[section]
\newtheorem{Thm}{Theorem}[section]
\newtheorem{Def}{Definition}[section]

\begin{document}
\begin{center}{\bf\Large Non-homogeneous fractional Schr\"odinger equation}\medskip

\bigskip

\bigskip

{C\'esar E. Torres Ledesma}

 Departamento de Matem\'aticas, \\
 Universidad Nacional de Trujillo,\\
 Av. Juan Pablo II s/n. Trujillo-Per\'u\\
 {\sl  (ctl\_576@yahoo.es, ctorres@dim.uchile.cl)}

\end{center}

\medskip

\medskip
\medskip
\medskip
\medskip

\begin{abstract}
In this article we are interested in the non-homogeneous fractional Schr\"odinger equation  
\begin{eqnarray}\label{eq00}
&(-\Delta)^{\alpha}u(x) + V(x)u(x) = f(u) + h(x)   \mbox{ in } \mathbb{R}^{n}. 
\end{eqnarray}
By using mountain pass Thoerem and Ekeland's variational principle, we prove the existence of two solutions for (\ref{eq00}).
\end{abstract}

\date{}

\setcounter{equation}{0}
\section{ Introduction}

Recently, a great attention has been focused on the study of problems involving the fractional Laplacian, from a pure mathematical point of view as well as from concrete applications, since this operator naturally arises in many different contexts, such as, obstacle problems, financial mathematics, phase transitions, anomalous diffusions, crystal dislocations, soft thin films, semipermeable membranes, flame propagations, conservation laws, ultra relativistic limits of quantum mechanics, quasi-geostrophic flows, minimal surfaces, materials science and water waves. The literature is too wide to attempt a reasonable list of references here, so we derive the reader to the work by Di Nezza, Patalluci and Valdinoci \cite{EDNGPEV}, where a more extensive bibliography and an introduction to the subject are given.

In the context of fractional quantum mechanics, non-linear fractional Schr\"odinger equation has been proposed by   Laskin \cite{NL-1}, \cite{NL-2} as a result of expanding the Feynman path integral, from the Brownian-like to the L\'evy-like quantum mechanical paths. In the last 10 years, there has been a lot of interest in the study of the fractional Schr\"odinger equation
\begin{equation}\label{I01}
(-\Delta)^{\alpha}u + V(x)u = f(x,u)\;\;\mbox{in} \;\;\mathbb{R}^{n}.
\end{equation}
where the nonlinearity $f$ satisfies some general conditions. See, for instance, Feng \cite{BF}, Chang \cite{XC-1}, \cite{XC-2}, Cheng \cite{MC}, Dipiero, Palatucci and Valdinoci \cite{SDGPEV}, Dong and Xu \cite{JDMX}, Felmer, Quaas and Tan \cite{PFAQJT}, de Oliveira, Costa and Vaz \cite{EOFCJV} and Secchi \cite{SS-1}, \cite{SS-2}. 

To the author's knowledge, most of these works assumed that there exists a trivial solution, namely $0$, for (\ref{I01}). There seems to have been very little progress on existence theory for (\ref{I01}) without trivial solutions.

This paper studies the existence of solutions $u\in H^{\alpha}(\mathbb{R}^{n})$ for the fractional equation
\begin{equation}\label{I02}
(-\Delta)^{\alpha}u(x) + V(x)u(x) = f(u) + h(x)   \mbox{ in } \mathbb{R}^{n},\;\;n\geq 2,
\end{equation}
where $0 < \alpha < 1$, $(-\Delta)^{\alpha}$ stands for the fractional laplacian defined by
$$
(-\Delta)^{\alpha}u(x) = p.v. \int_{\mathbb{R}^{n}} \frac{u(x) - u(z)}{|x-z|^{n+2\alpha}}dz.
$$
This problem is a model for $(\ref{I01})$ without trivial solutions and present specific mathematical difficulties. 

When $\alpha = 1$ in (\ref{I02}), we have the classical non homogeneous nonlinear  Schr\"odinger equation, which has been studied extensively by many authors in the last few decades, see, for example \cite{DCHZ}, \cite{ZWHZ}, \cite{XZ},  and references therein, where the existence and multiplicity results have been studied.

Troughout the paper we assume that
\begin{itemize}
\item[($V_{1}$)] $V\in C(\mathbb{R}^{n}, \mathbb{R})$ and there exists a constant $V_{0}>0$ such that $V(x) \geq V_{0},\;\;\forall x\in \mathbb{R}^{n}$,
\item[$(V_{2})$] $\lim_{|x|\to \infty} V(x) = \infty$,
\end{itemize}
Regarding $f$ we consider
\begin{itemize}
\item[$(f_{1})$] $f\in C(\mathbb{R})$, $f(0) = 0$,
\item[$(f_{2})$] $f(t) = o(|t|)$ as $t\to 0$,
\item[$(f_{3})$] $f(t) = o(|t|^{\frac{n+2\alpha}{n-2\alpha}})$ as $t \to \infty$,
\item[$(f_{4})$] There is a constant $\mu >2$ such that
\begin{equation}\label{I03}
0< \mu F(u) = \int_{0}^{u}f(s)ds \leq uf(u),\;\;u\not \equiv 0,
\end{equation}
\end{itemize}
and for $h$ we consider
\begin{itemize}
\item[(H)] $h\in L^{2}(\mathbb{R}^{n})$, $h \not\equiv 0$ and 
$$
\|h\|_{L^{2}(\mathbb{R}^{n})} < \frac{\varrho}{C_{e}} \left( \frac{1}{2} - \epsilon C_{e}^{2} + (\epsilon + K_{\epsilon}) C_{e}^{2_{\alpha}^{*}} \varrho ^{2_{\alpha}^{*}-2}\right),
$$
where $\varrho>0$ is given by the first geometrical condition of the mountain pass theorem and $C_{e}$ is the Sobolev constant.
\end{itemize}

 Our main result is as follows.
 \begin{Thm}\label{Ptm01}
 Under assumptions $(V_{1})-(V_{2})$, $(f_{1}) - (f_{4})$ and $(H)$, (\ref{I02}) hast at least two solutions.
 \end{Thm}

Our study is motivated by \cite{MC}, \cite{SS-1}, \cite{SS-2}. In \cite{MC} Cheng proved the existence of bound state solutions to (\ref{I01}) with $f(t) = t^{q}$ and unbounded potential by using Lagrange multiplier method and Nehari's manifold approach. It is worth noticing that under the assumption that potential $V(x) \to \infty$ as $|x| \to \infty$, the embedding $H_{V}^{\alpha}(\mathbb{R}^{n}) \hookrightarrow L^{q}(\mathbb{R}^{n})$ is compact, where
$$
H_{V}^{\alpha}(\mathbb{R}^{n}) = \left\{ u\in H^{\alpha}(\mathbb{R}^{n})/\;\; \int_{\mathbb{R}^{n}} V(x)u^{2}(x)dx < \infty\right\}
$$
and $2\leq q < \frac{2n}{n-2\alpha}$. In \cite{SS-1} Secchi has studied the equation (\ref{I01}). Under the same assumption on $V$, the existence of a ground states is obtained by Mountain pass Theorem. In \cite{SS-2}, Secchi looks for a radially symmetric solution of (\ref{I02}), with $f$ does not depend on $x$, namely de considered
$$
(-\Delta)^{\alpha}u(x) + V(x)u(x) = f(u)
$$
where the nonlinearity $f$ satisfies rather weak assumptions , which are comparable to those in \cite{HBPL}. By using the monoticity trick of Struwe-Jeanjean, Secchi shows the existence of radial solution.

Our theorem extends these result to the case $h\neq 0$. Under this assumption, the problem of existence of solutions is much more delicate, because the extra difficulties arise in studying the properties of the corresponding action functional $I: H_{V}^{\alpha}(\mathbb{R}^{n}) \to \mathbb{R}$

The problem here is as follows. We are given two sequence of almost critical point in $H_{V}^{\alpha}(\mathbb{R}^{n})$. The first one, obtained by Ekeland's variational principle, is contained in a small ball centered at $0$. Using the mountain pass geometry of the action functional, the existence of the second sequence is established. Both sequence are weakly convergent in $H_{V}^{\alpha}(\mathbb{R}^{n})$. The question is whether their limits are equal to each other or they define two geometrically distinct solutions of (\ref{I02}). The PS-condition is enough to obtain two solutions. The assumption ($V_{2}$) ensure the PS-condition at each level. In fact one needs the PS-condition only at two levels.

This article is organized as follows.  In Section $\S 2$ we present preliminaries with the main tools and the functional setting of the problem. In Section $\S 3$ we prove the Theorem \ref{Ptm01}.

\section{Preliminaries}

In this section, we collet some information to be used in the paper. Sobolev spaces of fractional order are the convenient setting for our equation. A very complete introduction to fractional Sobolev spaces can be found in \cite{EDNGPEV}.

We recall that the fractional Sobolev space $H^{\alpha}(\mathbb{R}^{n})$ is defined for any $\alpha \in (0,1)$ as
$$
H^{\alpha}(\mathbb{R}^{n}) = \left\{ u\in L^{2}(\mathbb{R}^{n})/\;\;\int_{\mathbb{R}^{n}}\int_{\mathbb{R}^{n}} \frac{|u(x) - u(z)|^{2}}{|x-z|^{n+2\alpha}} < \infty \right\}.
$$
This space is endowed with the Gagliardo norm
$$
\|u\|_{\alpha}^{2} = \int_{\mathbb{R}^{n}} u^{2}(x) dx + \int_{\mathbb{R}^{n}}\int_{\mathbb{R}^{n}} \frac{|u(x) - u(z)|^{2}}{|x-z|^{n+2\alpha}}.
$$

Regarding the space $H^\alpha(\mathbb{R}^n)$ we recall the following embedding theorem, whose proof can be found in \cite{EDNGPEV}.
\begin{Thm}\label{Ptm1}
Let $\alpha \in (0,1)$, then there exists a positive constant $C = C(n,\alpha)$ such that
\begin{equation}\label{P01}
\|u\|_{L^{2_{\alpha}^{*}}(\mathbb{R}^{n})}^{2} \leq C \int_{\mathbb{R}^{n}} \int_{\mathbb{R}^{n}} \frac{|u(x) - u(y)|^{2}}{|x-y|^{n+2\alpha}}dydx
\end{equation}
and then we have that  $H^{\alpha}(\mathbb{R}^{n}) \hookrightarrow L^{q}(\mathbb{R}^{n})$ is continuous for all $q \in [2, 2_{\alpha}^{*}]$.  

Moreover, $H^{\alpha}(\mathbb{R}^{n}) \hookrightarrow L^{q}(\Omega)$ is compact for any bounded set $\Omega\subset \mathbb{R}^{n}$ and for all $q \in [2, 2_{\alpha}^{*})$, where $2_{\alpha}^{*} = \frac{2n}{n-2\alpha}$ is the  critical exponent.
\end{Thm}

Now we consider the Hilbert space $H_{V}^{\alpha}(\mathbb{R}^{n})$ defined by
$$
H_{V}^{\alpha} = \left\{ u\in H^{\alpha}(\mathbb{R}^{n})/\;\; \int_{\mathbb{R}^{n}} V(x) u^{2}(x)dx < \infty \right\}
$$
endowed with the inner product
$$
\langle u, w \rangle_{V} = \int_{\mathbb{R}^{n}}\int_{\mathbb{R}^{n}} \frac{[u(x) - u(z)][w(x) - w(z)]}{|x-z|^{n+2\alpha}} + \int_{\mathbb{R}^{n}} V(x) u(x)w(x)dx,
$$ 
and norm
$$
\|u\|_{V}^{2} = \int_{\mathbb{R}^{n}}\int_{\mathbb{R}^{n}} \frac{|u(x) - u(z)|^{2}}{|x-z|^{n+2\alpha}} + \int_{\mathbb{R}^{n}} V(x)u^{2}(x)dx
$$

By $(V_{1})$ it is standard to prove that $H_{V}^{\alpha}(\mathbb{R}^{n})$ is continuously embedded in $H^{\alpha}(\mathbb{R}^{n})$ and by Theorem \ref{Ptm1}, we have that $H_{V}^{\alpha}(\mathbb{R}^{n}) \hookrightarrow L^{q}(\mathbb{R}^{n})$ is continuous for all $q\in [2, 2_{\alpha}^{*}]$. Moreover, we have the following compactness theorem
\begin{Thm}\label{Ptm2}
\cite{MC} Suppose that $(V_{1})$ and $(V_{2})$ holds. Then $H_{V}^{\alpha}(\mathbb{R}^{n}) \hookrightarrow L^{q}(\mathbb{R}^{n})$ is compact for all $q\in [2, 2_{\alpha}^{*})$.
\end{Thm}

Moreover we consider the following Lemma
\begin{Lem}\label{Plm1}
\cite{MGST} Suppose the $\beta >1$ and the function $f \in C(\mathbb{R})$ satisfies
$$
f(t) = o(|t|) \;\;\mbox{as}\;\;|t|\to 0\;\;\mbox{and}\;\;f(t) = o(|t|^{\beta}) \;\;\mbox{as}\;\;|t|\to \infty
$$
If $\{u_{k}\}_{k}$ is a bounded sequence in $L^{\beta + 1}(\mathbb{R}^{n})$ and $u_{k} \to u$ in $L^{2}(\mathbb{R}^{n})$ then
$$
\int_{\mathbb{R}^{n}} |f(u_{k})(u_{k} - u)|dx \to 0\;\;\mbox{as}\;\;k \to \infty.
$$
\end{Lem}
\section{Proof of theorem \ref{Ptm01}}

In this section, our goal is to prove the existence of solutions of equation (\ref{I02}). We start with a precise definition of the notion of solutions for equation (\ref{I02}).

\begin{Def}\label{def01}
We say that $u\in H_{V}^{\alpha}(\mathbb{R}^{n})$ is a weak solution of (\ref{I02}) if 
$$
\langle u, w\rangle_{V} = \int_{\mathbb{R}^{n}} ( f(u(x)) + h(x) )w(x)dx,\;\;\mbox{for all}\;\;v\in H_{V}^{\alpha}(\mathbb{R}^{n}).
$$
\end{Def}

We prove the existence of weak solution of (\ref{I02}) finding a critical point of the functional $I: H_{V}^{\alpha} (\mathbb{R}^{n}) \to \mathbb{R}$ defined by
\begin{equation}\label{eq01}
I(u) = \frac{1}{2}\|u\|_{V}^{2} - \int_{\mathbb{R}^{n}} F(u(x))dx - \int_{\mathbb{R}^{n}} h(x)u(x)dx.
\end{equation}
Using the properties of the Nemistky operators and the compact embedding Theorem \ref{Ptm2}, we can prove that the functional $I\in C^{1}(H_{V}^{\alpha}(\mathbb{R}^{n}), \mathbb{R})$ and we have
\begin{equation}\label{eq02}
I'(u)w = \langle u,w \rangle_{V} - \int_{\mathbb{R}^{n}} f(u(x))w(x)dx - \int_{\mathbb{R}^{n}}h(x)w(x)dx,\;\;\foral \;w\in H_{V}^{\alpha}(\mathbb{R}^{n})
\end{equation}

In order to prove Theorem \ref{Ptm01} we use the mountain pass Theorem (see \cite{PR} Theorem 2.2) and Ekeland's variational principle (see \cite{JMMW} Theorem 4.1 and Corollary 4.1). The proof will be divided into a sequence of Lemmas.
\begin{Lem}\label{lm01}
Suppose that $(V_{1})-(V_{2}), (f_{1})-(f_{4})$ and (H) holds. Then the functional $I: H_{V}^{\alpha}(\mathbb{R}^{n}) \to \mathbb{R}$ satisfies the Palais-Smale condition. 
\end{Lem}

\noindent 
{Proof.} Let $\{u_{k}\}$ be a sequence in $H_{V}^{\alpha}(\mathbb{R}^{n})$ such that
\begin{equation}\label{eq03}
|I(u_{k})| \leq C,\;\; I'(u_{k}) \to 0\;\;\mbox{in}\;\;(H_{V}^{\alpha}(\mathbb{R}^{n}))^{*}\;\;\mbox{as}\;\;k\to \infty.
\end{equation}
There exists $k_{0}$ such that for $k\geq k_{0}$
$$
|I'(u_{k})u_{k}| \leq \|u_{k}\|_{V}.
$$
Then
\begin{eqnarray*}
C + \|u_{k}\|_{V} & \geq & I(u_{k}) - \frac{1}{\mu} I'(u_{k})u_{k}\\
& = & \left(\frac{1}{2} - \frac{1}{\mu} \right)\|u_{k}\|_{V}^{2} +  \int_{\mathbb{R}^{n}} (\frac{1}{\mu}f(u_{k})u_{k} - F(u_{k}))dx \\
& & - \left(1- \frac{1}{\mu} \right) \int_{\mathbb{R}^{n}} h(x)u_{k}(x)dx\\
& \geq & \left( \frac{1}{2} - \frac{1}{\mu}\right)\|u\|_{V}^{2} - \left( C_{e} - \frac{C_{e}}{\mu} \right)\|h\|_{L^{2}}\|u_{k}\|_{V},
\end{eqnarray*}
so, $\{u_{k}\}$ is bounded in $H_{V}^{\alpha}(\mathbb{R}^{n})$. By Theorem \ref{Ptm2}, $H_{V}^{\alpha}(\mathbb{R}^{n}) \hookrightarrow L^{2}(\mathbb{R}^{n})$, compactly and $H_{V}^{\alpha}(\mathbb{R}^{n}) \hookrightarrow L^{2_{\alpha}^{*}}(\mathbb{R}^{n})$ continuously. Then, there exists a subsequence, still denoted by $\{u_{k}\}$ such that
\begin{eqnarray*}
&& u_{k} \rightharpoonup u \;\;\mbox{in}\;\;H_{V}^{\alpha}(\mathbb{R}^{n}),\\
&& u_{k} \to u \;\;\mbox{in}\;\;L^{2}(\mathbb{R}^{n}). 
\end{eqnarray*}  
By another hand
\begin{eqnarray*}
\|u_{k}\|_{V}^{2} - \langle u_{k}, u\rangle_{V} & = & I'(u_{k})(u_{k} - u) + \int_{\mathbb{R}^{n}} f(u_{k})(u_{k} - u)dx\\
&& + \int_{\mathbb{R}^{n}} f(x)(u_{k}(x) - u(x))dx.
\end{eqnarray*}
Hence, by Lemma \ref{Plm1} and (\ref{eq03})
$$
\lim_{k\to \infty} (\|u_{k}\|_{V}^{2} - \langle u_{k}, u\rangle_{V}) = 0.
$$
Then $\lim_{k\to \infty} \|u_{k}\|_{V} = \|u\|_{V}$ and, therefore, the sequence $\{u_{k}\}$ converges to $u$ strongly in $H_{V}^{\alpha}(\mathbb{R}^{n})$. $\Box$ 

\begin{Lem}\label{lm02}
Suppose that $(V_{1})-(V_{2}), (f_{1}) - (f_{4})$ and $(H)$ holds. There are $\varrho >0$ and $\tau>0$ such that
$$
I(u) \geq \tau\;\;\mbox{for}\;\;\|u\|_{V} = \varrho
$$ 
\end{Lem}

\noindent
{\bf Proof.} By continuous embedding
\begin{equation}\label{eq04}
\|u\|_{L^{2}(\mathbb{R}^{n})} \leq C_{e} \|u\|_{V},\;\;\|u\|_{L^{2_{\alpha}^{*}}(\mathbb{R}^{n})} \leq C_{e}\|u\|_{V}.
\end{equation}
By $(f_{2})$ and $(f_{3})$, for every $\epsilon$ there exists $\rho, \delta >0$ such that
\begin{eqnarray*}
&&|f(t)| \leq \epsilon |t|^{\frac{n+2\alpha}{n-2\alpha}}\;\;\mbox{for}\;\;|t|\geq \rho,\;\;\mbox{and}\\
&& |f(t)| \leq \epsilon |t|\;\;\mbox{for}\;\;|t| \leq \delta.
\end{eqnarray*} 
Therefore we have
$$
|f(t)| \leq \epsilon (|t| + |t|^{\frac{n+2\alpha}{n-2\alpha}}) + K_{\epsilon} |t|^{\frac{n+2\alpha}{n-2\alpha}},
$$
where $K_{\epsilon} = \delta^{-\frac{n+2\alpha}{n-2\alpha}} \max_{\delta \leq |t| \leq \rho}|f(t)|$. Then we have
\begin{equation}\label{eq05}
|F(t)| \leq \epsilon  ( |t|^{2} + |t|^{2_{\alpha}^{*}}) + K_{\epsilon}|t|^{2_{\alpha}^{*}}.
\end{equation}
Let $0< \epsilon < \frac{1}{2C_{e}^{2}}$. By (\ref{eq04}) and (\ref{eq05}) we have
\begin{eqnarray*}
\int_{\mathbb{R}^{n}} F(u)dx & \leq & \epsilon \left( \|u\|_{L^{2}(\mathbb{R}^{n})}^{2} + \|u\|_{L^{2_{\alpha}^{*}}(\mathbb{R}^{n})}^{2_{\alpha}^{*}} \right) + K_{\epsilon} \|u\|_{L^{2_{\alpha}^{*}}(\mathbb{R}^{n})}^{2_{\alpha}^{*}}\\
& \leq & \epsilon C_{e}^{2} \|u\|_{V}^{2} + (\epsilon + K_{\epsilon})C_{e}^{2_{\alpha}^{*}} \|u\|_{V}^{2_{\alpha}^{*}}
\end{eqnarray*}
and 
\begin{equation}\label{eq06}
I(u) \geq \left( \frac{1}{2} - \epsilon C_{e}^{2}\right)\|u\|_{V}^{2} - (\epsilon + K_{\epsilon})C_{e}^{2_{\alpha}^{*}}\|u\|_{V}^{2_{\alpha}^{*}} - C_{e}\|h\|_{L^{2}(\mathbb{R}^{n})} \|u\|_{V}.
\end{equation}
Taking $\|u\|_{V} = \varrho$ then $I(u) \geq \tau>0$ by (H). $\Box$

\begin{Lem}\label{lm03}
Suppose that $(V_{1})-(V_{2}), (f_{1}) - (f_{4})$ and $(H)$ holds. There is $e \in \overline{B(0, \varrho)}^{c}$ such that $I(e) \leq 0$
\end{Lem}

\noindent
{\bf Proof.} Since $h \not\equiv 0$, we can choose a function $\varphi \in H^{\alpha}(\mathbb{R}^{n})$ such that
$$
\int_{\mathbb{R}^{n}} h(x)\varphi (x)dx >0.
$$
By ($f_{4}$) it follows that there exists a constant $m>0$ such that
\begin{equation}\label{I04}
F(u) \geq m|u|^{\mu}\;\;\mbox{if}\;\;|u|\geq 1,
\end{equation}
so, for $\lambda \in (0, + \infty)$, we have
\begin{eqnarray*}
I(\lambda \varphi) & = & \frac{\lambda^{2}}{2} \|\varphi\|_{V}^{2} - \int_{\mathbb{R}^{n}} F(\lambda \varphi) dx - \lambda\int_{\mathbb{R}^{n}}h(x)\varphi (x)dx \\
& \leq & \frac{\lambda^{2}}{2} \|\varphi\|_{V}^{2} - m\lambda^{\mu}\int_{\{|u|\geq 1\}}|\varphi|^{\mu}dx - \lambda \int_{\mathbb{R}^{n}} h(x)\varphi (x)dx.
\end{eqnarray*}
Since $\mu > 2$, $I(\lambda \varphi) \to -\infty$ as $\lambda \to +\infty$. Hence, there is $\lambda \in (0, +\infty)$ such that
$$
\|\lambda \varphi\|_{V} > \varrho \;\;\mbox{and}\;\; I(\lambda \varphi) \leq 0.
$$ 
$\Box$

\noindent 
{\bf Proof of Theorem \ref{Ptm01}}

Since $I(0) = 0$ and $I$ satisfies Lemmas \ref{lm01} - \ref{lm03}, it follows by the mountain pass Theorem that $I$ has a critical value $c$ given by
$$
c = \inf_{\gamma \in \Gamma} \max_{t\in [0,1]} I(\gamma (t)),
$$
where $\Gamma = \{\gamma \in C([0,1], H_{V}^{\alpha}(\mathbb{R}^{n}))/\;\; \gamma (0) = 0, I(\gamma (1)) \leq 0\}$.  By definition , it follows that $c\geq \varrho >0$.

From (\ref{eq06}), we conclude that $I$ is bounded from below on $\overline{B(0, \varrho)}$. Set
\begin{equation}\label{eq07}
\overline{c} = \inf_{\|u\|_{V} \leq \varrho} I(u).
\end{equation}
Hence $I(0) = 0$ implies $\overline{c} \leq 0$. Thus $\overline{c} < c$. By Ekeland's variational principle, there is a minimizing sequence $\{w_{k}\} \subset \overline{B(0, \varrho)}$ such that
$$
I(w_{k}) \to \overline{c}\;\;\mbox{and}\;\; I'(w_{k}) \to 0\;\;\mbox{as} \;\;k\to \infty.
$$
From Lemma \ref{lm01}, $\overline{c}$ is a critical value of $I$. Consequently, $I$ has at least two critical points. $\Box$

\noindent {\bf Acknowledgements:}
C. T. was  partially supported by MECESUP 0607 and CMM.


\begin{thebibliography}{99}
\bibitem{AAPR}A. Ambrosetti and P. Rabinowitz, \emph{``Dual variational methods in critical points theory and applications''}, J. Func. Anal. {\bf 14}, 349-381(1973).

\bibitem{HBPL}H. Berestycki and P.-L. Lions, \emph{``Nonlinear scalar field equations. I. Existence of a ground state''}, Arch. Rational Mech. Anal. {\bf 82}, nº 4, 313-345(1983).


\bibitem{DCHZ}D. Cao and H. Zhou, \emph{``Multiple positive solutions of nonhomogeneous semilinear elliptic equation in $\mathbb{R}^{n}$''}, Proc. Roy. Soc. Edinburgh Sect. A, {\bf 126}, 443-463{1996}.

\bibitem{XC-1}X. Chang, \emph{``Ground state of fractional Schr\"odinger equation on $\mathbb{R}^{n}$''}, Proc. Edinb. Math. Soc. (to be published).

\bibitem{XC-2}X. Chang, \emph{``Ground state solutions of asymptotically linear fractional Schr\"odinger equations''}, J. Math. Phys., {\bf 54}, 061504(2013).

\bibitem{MC}M. Cheng, \emph{``Bound state for the fractional Schr\"odinger equation with undounded potential''}, J. Math. Phys. {\bf 53}, 043507 (2012).

\bibitem{DC}D. Costa, \emph{``An invitation to variational methods in differential equations''}, Birkhäuser, Boston, 2007.


\bibitem{EDNGPEV}E. Di Nezza, G. Patalluci and E. Valdinoci, \emph{``Hitchhiker's guide to the fractional Sobolev spaces''}, Bull. Sci. math., 2012.

\bibitem{SDGPEV}S. Dipierro, G. Palatucci and E. Valdinoci, \emph{``Existence and symmetry results for a Schr\"odinger type problem involving the fractional Laplacian''},  Le Matematiche, {\bf LXIII}, I, 201-216(2013).

\bibitem{JDMX}J. Dong and M.Xu, \emph{``Some solutions to the space fractional Schr\"odinger equation using momentum representation method''},
J. Math. Phys. {\bf 48}, 072105 (2007).

\bibitem{PFAQJT}P. Felmer, A. Quaas and J. Tan \emph{``Positive solutions of nonlinear Schr\"dinguer equation with the fractional laplacian''}, Proceedings of the Royal Society of Edinburgh: Section A Mathematics, {\bf 142}, No 6, 1237-1262(2012).

\bibitem{BF}B. Feng, \emph{``Ground states for the fractional Schr\"odinger equation''}, EJDE, No 127, 1-11(2013)

\bibitem{guan1} Q-Y. Guan  \emph{``Integration by Parts Formula for Regional Fractional Laplacian."} Commun. Math. Phys. 266, 289Ð329 (2006). 

\bibitem{guan2} Q-Y. Guan, Z.M. Ma  \emph{``The reflected $\alpha$-symmetric stable processes and regional fractional Laplacian.''} Probab. Theory Relat. Fields 134(4), 649Ð694 (2006)

\bibitem{XGMX}X. Guo and M. Xu, \emph{``Some physical applications of fractional Schr\"odinger equation''}, J. Math. Phys. {\bf 47}, 082104 (2006).

\bibitem{MGST}M. Grossinho and S. Tersian, \emph{``An introduction to minimax theorems and their applications to differential equations''}, Kluwer Academic Publishers, 2001. 

\bibitem{HIGN}H. Ishii and G. Nakamura, \emph{``A class of integral equations and approximation of p-Laplace equations''}, Calc. Var. {\bf 37}, 485-522(2010).

\bibitem{NL-1}N. Laskin, \emph{``Fractional quantum mechanics and L\'evy path integrals''}, Phys. Lett. A {\bf 268}, 298 - 305 (2000).

\bibitem{NL-2}N. Laskin, \emph{``Fractional Schr\"odinger equation''}. Phys. Rev. E {\bf 66}, 056108 (2002).

\bibitem{JMMW}J. Mawhin and M. Willen, \emph{``Critical point theory and Hamiltonian systems''}, Applied Mathematical Sciences 74, Springer, Berlin, 1989.


\bibitem{EOFCJV}E. de Oliveira, F. Costa, and J. Vaz, \emph{``The fractional Schr\"odinger equation for delta potentials''}, J. Math. Phys. {\bf 51},
123517 (2010).

\bibitem{PR}P. Rabinowitz, \emph{``Minimax method in critical point theory with applications to differential equations''}, CBMS Amer. Math. Soc., No {\bf 65}, 1986.

\bibitem{SS-1}S. Secchi, \emph{``Ground state solutions for nonlinear fractional Schrödinger equation in $\mathbb{R}^{n}$''},  J. Math. Phys., {\bf 54}, 031501(2013)

\bibitem{SS-2}S. Secchi, \emph{``On fractional Schrödinger equation in $\mathbb{R}^{n}$ without the Ambrosetti-Rabinowitz condition''}, preprint. arXiv:1210.0755v1 [math.AP] 2 Oct 2012.

\bibitem{ZWHZ}Z. Wang and H. Zhou, \emph{``Positive solutions for a nonhomogeneous elliptic equation on $\mathbb{R}^{n}$ without (AR) condition''}, J. Math. Anal. Appl., {\bf 353}, 470-479{2009}.

\bibitem{XZ}X. Zhu, \emph{``A perturbation result on positive entire solutions of a semilinear elliptic equation''}, J. Diff. Equ., {\bf 92}, 2, 162-178{1991}.



















\end{thebibliography}
\end{document}